\documentclass{article}   

\usepackage{amsmath}
\usepackage{amssymb,amsfonts}
\usepackage{graphicx}
\usepackage{times,amssymb,amscd}
\usepackage[USenglish]{babel}
\usepackage{verbatim}
\usepackage{enumerate}

\begin{document}

\title{Optimally Dense Packings for Fully Asymptotic Coxeter Tilings by Horoballs of Different Types
}

\author{Robert Thijs Kozma\\
		Department of Mathematics and Statistics\\
		Boston University, 111 Cummington St.\\ Boston, MA 02215 USA\\
		Present Address: Department of Mathematics\\
		SUNY Stony Brook\\
		Stony Brook, NY 11794-3651 USA\\
		email: rkozma@math.sunysb.edu\\
		\\
		\\
Jen\H{o}  Szirmai\\
Budapest University of Technology and Economics\\
		Institute of Mathematics, Department of Geometry \\
		H-1521 Budapest, Hungary \\
		email: szirmai@math.bme.hu
}

\maketitle

\begin{abstract}
The goal of this paper is to determine the optimal horoball packing arrangements and their densities
for all four fully asymptotic Coxeter tilings
(Coxeter honeycombs) in hyperbolic 3-space $\mathbb{H}^3$. 
Centers of horoballs are required to lie at vertices of the
regular polyhedral cells constituting the tiling. We allow horoballs of
different types at the various vertices. Our results are derived through
a generalization of the projective methodology for hyperbolic spaces. The
main result states that  the known B\"or\"oczky--Florian 
density upper bound for ``congruent horoball" packings of $\mathbb{H}^3$
remains valid for the class of fully asymptotic Coxeter tilings, even if
packing conditions are relaxed by allowing for horoballs of different
types under prescribed symmetry groups. The consequences of this remarkable result are discussed for
various Coxeter tilings.
\end{abstract}

\newtheorem{theorem}{Theorem}[section]
\newtheorem{corollary}[theorem]{Corollary}
\newtheorem{lemma}[theorem]{Lemma}
\newtheorem{exmple}[theorem]{Example}
\newtheorem{defn}[theorem]{Definition}
\newtheorem{rmrk}[theorem]{Remark}
\newtheorem{proposition}[theorem]{Proposition}
\newenvironment{definition}{\begin{defn}\normalfont}{\end{defn}}
\newenvironment{remark}{\begin{rmrk}\normalfont}{\end{rmrk}}
\newenvironment{example}{\begin{exmple}\normalfont}{\end{exmple}}
\newenvironment{acknowledgement}{Acknowledgement}


\section{Introduction}

Local optimal ball packings for regular tilings have been studied extensively in the literature. 
Of special interest are tilings in hyperbolic
$n$-space; B\"or\"oczky and Florian \cite{B--F64} gave the
universal density upper bound for all congruent ball packings of the 3-dimensional hyperbolic space $\mathbb{H}^3$ without any symmetry
assumptions. This classical result provides the density upper bound realized by a regular horoball packing of (3, 3, 6) in $\mathbb{H}^3$, 
as shown in Section 2. The
optimal density is related to the Dirichlet-Voronoi cell of every
ball, as follows:
$$
s_0=( 1+ \frac{1}{2^2}-\frac{1}{4^2}-\frac{1}{5^2}+\frac{1}{7^2}+\frac{1}{8^2}- \ -
\ + \ + \dots)^{-1} \approx 0.85327609.
$$
This limit is achieved by 4 horoballs  centered at the vertices of a regular
ideal simplex, tangent to each other at the ``midpoints" of the edges, i.e. as the projection of the simplex center
into any edge of the simplex.

Papers \cite{Be}, \cite{Bo--R}, \cite{G--K--K}, \cite{K91},
\cite{R06} introduce recent novel developments in the
classical topic of the ball (or sphere) packings of $\mathbb{H}^3$.
Locally dense (optimal) ball, horoball packings in $\mathbb{H}^3$
are of great significance, as important information 
regarding crystal structures can be obtained using locally optimal ball and
horoball arrangements.

The present work is based on the projective interpretation of the
hyperbolic geometry, proposed in \cite{M93}, \cite{M97}. In subsequent works the second author studied a class of face
transitive tilings  \cite{D--H--M}, the so-called generalized Lambert-cube tilings in 
\cite{Sz03-1} and \cite{Sz05-1}, where an algorithmic approach for
determining volumes of hyperbolic polyhedra was developed and
implemented. Using this novel approach, the locally optimal ball
packings were found for the configurations in which the ball and
horoball centers lie either within the Lambert-cubes or at the
vertices of the cubes, respectively, and the optimal packing densities of the
corresponding tilings were computed \cite{Sz05-1}. Optimal ball
and horoball packings of the regular Coxeter honeycombs in
$\mathbb{H}^d, (d \ge 3)$  with one horoball type were found in 
\cite{Sz05-2}, \cite{Sz07-1}, and the optimal ball and horoball
packing densities computed.

In related work,  $d$-dimensional $(d\ge 3)$
hyperbolic prism honeycombs generated by ``inscribed hyperspheres" were investigated in
\cite{Sz06-1} and \cite{Sz06-2}. The optimal hyperball packings of
infinitely many 3-dimensional prism tilings (mosaics) together with their
metric data were determined in \cite{Sz06-2}. In hyperbolic 4-space
$\mathbb{H}^4$ there are only 2 honeycombs with metric data
corresponding to their 3-dimensional counterparts. The densities of the
optimal hyperball packings in 4-space are determined in \cite{Sz06-2}. In
$\mathbb{H}^5$ there are 3 types of such mosaics, and the
corresponding problems are extensively studied. In $\mathbb{H}^d \
(d>5)$ there are no longer any regular prism tilings.

In this paper, we study locally optimal ball and horoball packings in
the four fully asymptotic hyperbolic tilings of $\mathbb{H}^3$, while allowing different
types of horoballs to be centered at the vertices of the honeycombs. In
Section 2, we provide preliminaries on the $d$-dimensional honeycombs.
In Section 3, we introduce the projective model \cite{Sz05-1} to determine
the densities of the optimally dense horoball packings in hyperbolic
space $\mathbb{H}^d$. In Section 4 we determine the
optimal packing densities in $\mathbb{H}^3$ for various honeycombs,
when horoballs of various types are allowed. We find that the
densest possible packings yield density values identical to that of the
B\"or\"oczky--Florian bound \cite{B--F64}. In all
studied configurations, the optimal densities never surpass the B\"or\"oczky--Florian
upper bound, even when replacing the ``congruency" constraints with regularity constraints. We finish
the paper with conclusions and directions for future research.

\section{Overview on $d$-dimensional hyperbolic honeycombs}

Hyperbolic geometry is based on the principles of
Bolyai-Lobachevsky geometry \cite{P06}.  
A $d$-dimensional honeycomb $\mathcal{P}$, also referred to as a
solid tessellation or tiling, is an infinite collection of congruent
polyhedra (polytopes) that fit together face-to-face to fill the entire geometric
space $(\mathbb{H}^d~ (d \geqq 2))$ exactly once. 
We take the cells to be
congruent regular polyhedra. A honeycomb with cells congruent to a
given regular polyhedron $P$ exists if and only if the dihedral
angle of $P$ is a submultiple of $2\pi$ (in the hyperbolic plane
zero angles are also permissible). A complete classification of
honeycombs with bounded cells was first given by {\sc{Schlegel}} in
$1883$. The classification was completed by including the polyhedra
with unbounded cells, namely the fully asymptotic ones by
{\sc{Coxeter}} in 1954 \cite{C56}. {\it Such honeycombs exist only
for $d \le 5$} in hyperbolic $d$-space $\mathbb{H}^d$. In this paper Coxeter
honeycombs or Coxeter tilings refer to tilings described in Table 1.

An alternative approach to describing honeycombs involves analysis
of their symmetry groups. If $\mathcal{P}$ is a Coxeter honeycomb,
then any rigid motion moving one cell into another maps the entire
honeycomb onto itself. The symmetry group of a honeycomb is denoted
by $Sym \mathcal{P}$. The characteristic simplex $\mathcal{F}$ of
any cell $P \in \mathcal{P}$ is a fundamental domain of the symmetry
group $Sym \mathcal{P}$ generated by reflections in its facets which
are $(d-1)$-dimensional hyperfaces.

The scheme of a regular polytope $P$ is a weighted graph (diagram)  
characterizing $P \subset \mathbb{H}^d$ up to congruence. The nodes
of the scheme, numbered by $0,1,\dots,d$, correspond to the bounding
hyperplanes of $\mathcal{F}$. Two nodes are joined by an edge if the
corresponding hyperplanes are non-orthogonal. Let the set of weights
$(n_1,n_2,$ $n_3,\dots,n_{d-1})$ be the Schl\"{a}fli symbol of $P$,
and $n_d$ be the weight describing the dihedral angle of $P$, such
that the dihedral angle is equal to $ \frac{2\pi}{n_d}$. In this
case $\mathcal{F}$ is the Coxeter simplex with the scheme:

\begin{figure}[ht]
\centering
\includegraphics[width=7cm]{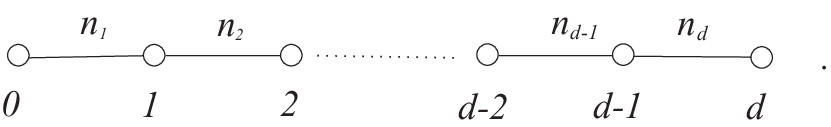}
\caption{Coxeter-Schl\"afli simplex scheme}
\end{figure}

The Schl\"{a}fli symbol of the honeycomb $\mathcal{P}$ is the
ordered set $(n_1,n_2,n_3,\dots,$ $n_{d-1},n_d)$ above.  A $(d+1)\times(d+1)$
symmetric matrix $(b^{ij})$ is constructed for each scheme in the
following manner: $b^{ii}=1$ and if $i \ne j\in \{0,1,2,\dots,d \}$
then $b^{ij}=-\cos{\frac{\pi}{n_{ij}}}$. For all angles between the
facets $i$,$j$ of $\mathcal{F}$ holds then $n_k= n_{k-1,k}$. Reversing the
numbering of the nodes of scheme $\mathcal{P}$ while keeping the
weights, leads to the scheme of the dual honeycomb $\mathcal{P}^*$
whose symmetry group coincides with $Sym \mathcal{P}$.

In this paper we investigate regular Coxeter honeycombs and their
optimal horoball packings in the hyperbolic space $\mathbb{H}^3$,
where the horoballs are allowed to be of different types. $Sym
\mathcal{P}$ denotes the symmetry group of the honeycomb
$\mathcal{P}_{n_1 n_2 \dots n_d}$, thus
$$
P_{n_1 n_2 \dots n_d}=\{\bigcup_{\gamma ~ \in ~ Sym \mathcal{P}_{n_1 n_2 \dots n_{d-1}}}
\gamma(\mathcal{F}_{n_1 n_2 \dots n_d}) \}.
$$
In order to calculate the packing density, we relate each ball or
horoball, respectively, to its regular polytope $P_{n_1 n_2 \dots
n_d}$ in which it is contained. These polytopes are not necessarily
assumed to be Dirichlet-Voronoi cells.
\begin{center}
Table 1: Classification of 3-dimensional Coxeter tilings

\vspace{3mm} \begin{tabular}{|l|l|l|}
\hline
No.& Description & Schl\"{a}fli symbol $(p, \ q, \ r)$  \\
\hline
1. & Cells having proper centers and vertices &  (3,5,3), \ (4,3,5), \ (5,3,4), \ (5,3,5) \\
\hline
2. & Fully asymptotic cells &  (3,3,6), \ (3,4,4), \ (4,3,6), \ (5,3,6) \\
\hline
3. & Infinite centers and proper  &  (3,6,3), \ (4,4,4), \ (6,3,6)\\
 & or ideal vertices &  (4,4,3), \ (6,3,3), \ (6,3,4), \ (6,3,5) \\
\hline
\end{tabular}
\end{center}

As listed in Table 1, Coxeter tilings with parameters in row 1
include cells having proper centers and vertices. The polyhedra of
honeycombs with Schl\"{a}fli symbols in row 3 of Table 1, have
infinite centers and proper or ideal vertices. The polyhedra of
tilings in row 2 of Table 1 are called fully asymptotic; moreover
their centers are proper and their vertices lie on the absolute
of the hyperbolic space i.e. they are ideal vertices.

\section{The projective model}

\subsection{Basic Notions}

Let X denote one of either the $d$-dimensional sphere
$\mathbb{S}^d$, the $d$ dimensional Euclidean space $\mathbb{E}^d$,
or the hyperbolic space $\mathbb{H}^d$, $d \geq 2$. For
$\mathbb{H}^d$ we use the projective model in Lorentz space
$\mathbb{E}^{1,d}$ of signature $(1,d)$, i.e.~$\mathbb{E}^{1,d}$ is
the real vector space $\mathbf{V}^{d+1}$ equipped with the bilinear
form of signature $(1,d)$
\begin{equation}
\langle ~ \mathbf{x},~\mathbf{y} \rangle = -x^0y^0+x^1y^1+ \dots + x^d y^d \tag{3.1}
\end{equation}
where the non-zero vectors
$$
\mathbf{x}=(x^0,x^1,\dots,x^d)\in\mathbf{V}^{d+1} \ \  \text{and} \ \ \mathbf{y}=(y^0,y^1,\dots,y^d)\in\mathbf{V}^{d+1},
$$
are determined up to real factors and they represent points in
$\mathcal{P}^d(\mathbb{R})$. $\mathbb{H}^d$ is represented as the
interior of the absolute quadratic form
\begin{equation}
Q=\{[\mathbf{x}]\in\mathcal{P}^d | \langle ~ \mathbf{x},~\mathbf{x} \rangle =0 \}=\partial \mathbb{H}^d \tag{3.2}
\end{equation}
in real projective space $\mathcal{P}^d(\mathbf{V}^{d+1},
\mbox{\boldmath$V$}\!_{d+1})$. All proper interior points $\mathbf{x} \in \mathbb{H}^d$ are characterized by
$\langle ~ \mathbf{x},~\mathbf{x} \rangle < 0$.

The points on the boundary $\partial \mathbb{H}^d $ in
$\mathcal{P}^d$ represent the absolute points at infinity of $\mathbb{H}^d $.
Points $\mathbf{y}$ with $\langle ~ \mathbf{y},~\mathbf{y} \rangle >
0$ lie outside $\partial \mathbb{H}^d $ and are called outer points
of $\mathbb{H}^d $. Let $P([\mathbf{x}]) \in \mathcal{P}^d$; a point
$[\mathbf{y}] \in \mathcal{P}^d$ is said to be conjugate to
$[\mathbf{x}]$ relative to $Q$ when $\langle ~
\mathbf{x},~\mathbf{y} \rangle =0$. The set of all points conjugate
to $P([\mathbf{x}])$ form a projective (polar) hyperplane
\begin{equation}
pol(P):=\{[\mathbf{y}]\in\mathcal{P}^d | \langle ~ \mathbf{x},~\mathbf{y} \rangle =0 \}. \tag{3.3}
\end{equation}
Hence the bilinear form $Q$ by (3.1) induces a bijection
(linear polarity $\mathbf{V}^{d+1} \rightarrow
\mbox{\boldmath$V$}\!_{d+1})$)
from the points of $\mathcal{P}^d$
onto its hyperplanes.

Point $X [\bold{x}]$ and the hyperplane $\alpha
[\mbox{\boldmath$a$}]$ are called incident if the value of the
linear form $\mbox{\boldmath$a$}$ on the vector $\bold{x}$ is equal
to zero; i.e., $\bold{x}\mbox{\boldmath$a$}=0$ ($\mathbf{x} \in \
\mathbf{V}^{d+1} \setminus \{\mathbf{0}\}, \ \mbox{\boldmath$a$} \in
\mbox{\boldmath$V$}_{d+1} \setminus \{\mbox{\boldmath$0$}\}$).
Straight lines in $\mathcal{P}^d$ are characterized by the
2-subspaces of $\mathbf{V}^{d+1}$ or $(d-1)$-spaces of $\
\mbox{\boldmath$V$}\!_{d+1}$ \cite{M97}.

Let $P \subset \mathbb{H}^d $ denote a polyhedron bounded by
hyperplanes $H^i$, which are characterized by unit normal vectors
$\mbox{\boldmath$b$}^i \in \mbox{\boldmath$V$}\!_{d+1}$ directed
inwards with respect to $P$:
\begin{equation}
H^i:=\{\mathbf{x} \in \mathbb{H}^d | \langle ~ \mathbf{x},~\mbox{\boldmath$b$}^i \rangle =0 \} \ \ \text{with} \ \
\langle \mbox{\boldmath$b$}^i,\mbox{\boldmath$b$}^i \rangle = 1. \tag{3.4}
\end{equation}
{\it We always assume $P$ to be an acute-angled polyhedron and the vertices to be proper points or to lie at infinity.}

The Gram matrix $G(P):=( \langle \mbox{\boldmath$b$}^i,
\mbox{\boldmath$b$}^j \rangle ) ~ {i,j \in \{ 0,1,2 \dots d \} }$ of
normal vectors $\mbox{\boldmath$b$}^i$ associated with $P$ is an
indecomposable symmetric matrix of signature $(1,d)$ with entries
$\langle \mbox{\boldmath$b$}^i,\mbox{\boldmath$b$}^i \rangle = 1$
and $\langle \mbox{\boldmath$b$}^i,\mbox{\boldmath$b$}^j \rangle
\leq 0$ for $i \ne j$, having the following geometrical meaning
\small
$$
\langle \mbox{\boldmath$b$}^i,\mbox{\boldmath$b$}^j \rangle =
\left\{
\begin{aligned}
&0 & &\text{if}~H^i \perp H^j,\\
&-\cos{\alpha^{ij}} & &\text{if}~H^i,H^j ~ \text{intersect \ on $P$ \ at \ angle} \ \alpha^{ij}, \\
&-1 & &\text{if}~\ H^i,H^j ~ \text{are parallel in hyperbolic sense}, \\
&-\cosh{l^{ij}} & &\text{if}~H^i,H^j ~ \text{admit a common perpendicular of length} \ l^{ij}.
\end{aligned}
\right.
$$
\normalsize

\begin{definition}
\cite{B--H}, \cite{K91} An orthoscheme $\mathcal{O}$ in X is a simplex
bounded by $d+1$ hyperplanes $H^0,\dots,H^d$ such that

$$
H^i \bot H^j, \  \text{for} \ j\ne i-1,i,i+1.
$$
\end{definition}

A plane orthoscheme is a right-angled triangle, the area of which
can be expressed by the defect formula. For orthoschemes we denote
the $(d-1)$-hyperface opposite to the vertex $A_i$ by $H^i$ $(0 \le
i \le d)$. An orthoscheme $\mathcal{O}$ has $d$ dihedral angles 
different from right angles. Let $\alpha^{ij}$ denote the dihedral
angle of $\mathcal{O}$ between the faces $H^i$ and $H^j$. Then we
have
\begin{equation}
\alpha^{ij}=\frac{\pi}{2}, \ \ \text{if} \ \ 0 \le i < j -1 \le d. \notag
\end{equation}

The remaining $d$ dihedral angles $\alpha^{i,i+1}, \ (0 \le i \le
d-1)$ are called the essential angles of $\mathcal{O}$. The initial
vertex $A_0$ and final vertex $A_d$ of the orthogonal edge-path
$$
\bigcup_{i=0}^{d-1} A_iA_{i+1}
$$
are called principal vertices of the orthoscheme.

In this work, the characteristic simplex $\mathcal{F}$ of any
honeycomb $\mathcal{P}$ with Schl\"afli symbol
$(n_1,n_2,n_3,\dots,n_{d})$ is an orthoscheme.

The matrix $(b^{ij})=G(P)$ is the so called Coxeter-Schl\"afli
matrix of the orthoscheme $\mathcal{F}$ with parameters
$n_1,n_2,n_3,\dots,n_{d}$:
\[
(b^{ij}):=\begin{pmatrix}
1& -\cos{\frac{\pi}{n_1}}& 0 &\dots&0 \\
-\cos{\frac{\pi}{n_1}} & 1 &-\cos{\frac{\pi}{n_2}}&\dots&0 \\
0 & -\cos{\frac{\pi}{n_2}} & 1 &\dots&0 \\
0 & 0 & -\cos{\frac{\pi}{n_3}} & \dots&0 \\
\hdotsfor[1.8]{5} \\
0 & \dots& 0 & -\cos{\frac{\pi}{n_d}} & 1
\end{pmatrix}. \tag{3.5}
\]

Inverting the Coxeter-Schl\"afli matrix $(b^{ij})$ (3.5) of an
orthoscheme gives the matrix $(a_{ij})$, which can be used to
express distances between two vertices through the formula
\cite{Sz06-2}:

\[
\cosh{\frac{d_{ij}}{k}}=\frac{-a_{ij}}{\sqrt{a_{ii} a_{jj}}}.   \tag{3.6}
\]
In this paper we set the sectional curvature of $\mathbb{H}^d$,
$K=-k^2$, to be $k=1$. The distance $s$ of two proper points
$(\mathbf{x})$ and $(\mathbf{y})$ is calculated by the formula:
\begin{equation}
\cosh{\frac{s}{k}}=\frac{-\langle ~ \mathbf{x},~\mathbf{y} \rangle }{\sqrt{\langle ~ \mathbf{x},~\mathbf{x} \rangle
\langle ~ \mathbf{y},~\mathbf{y} \rangle }} . \tag{3.7}
\end{equation}


\subsection{Characterization
of horoballs in Hyperbolic Space $\mathbb{H}^3$}

\begin{definition}
A horosphere in the hyperbolic geometry is the surface orthogonal to
the set of parallel lines, passing through the same point on the absolute quadratic surface 
(simply absolute) of  the hyperbolic space.
\end{definition}

We represent hyperbolic space $\mathbb{H}^3$ in the Cayley-Klein
ball model. We introduce a projective coordinate system using vector
basis $\bold{b}_i \ (i=0,1,2,3)$ for $\mathcal{P}^3$ where the
coordinates of center of the model is $A_2(1,0,0,0)$. We pick an
arbitrary point at infinity to be $A_3(1,0,0,1)$.

As it is known, the equation of a horosphere with center
$A_3(1,0,0,1)$ through point $S(1,0,0,s)$ is derived using the
surface pencil of the absolute sphere and a plane tangent to the sphere at
point $A_3(1,0,0,1)$. The equation of the absolute sphere is
$-x^0 x^0 +x^1 x^1+x^2 x^2+x^3 x^3=0.$ The equation of a plane tangent to the
absolute of our model at point $A_3(1,0,0,1)$ is
$x^0-x^3=0$.

The general equation of the horosphere is
\begin{equation}
0=\lambda (-x^0 x^0 +x^1 x^1+x^2 x^2+x^3 x^3)+\mu{(x^0-x^3)}^2. \notag
\end{equation}
This passes through point $S(1,0,0,s)$ so we may write
\begin{equation}
\lambda (-1+s^2)+\mu {(-1+s)}^2=0 \Rightarrow \frac{\lambda}{\mu}=\frac{{(s-1)}^2}{1-s^2} \notag
\end{equation}
\begin{gather}
\text{If $s \neq \pm1$, then} \ \ \  \frac{{(s-1)}^2}{1-s^2}(-x^0 x^0 +x^1 x^1+x^2 x^2+x^3 x^3)+{(x^0-x^3)}^2=0 \Leftrightarrow \notag \\
\Leftrightarrow (s-1)(-x^0 x^0 +x^1 x^1+x^2 x^2+x^3 x^3)-(1+s){(x^0-x^3)}^2=0 \notag
\end{gather}

This way we obtain the following equation for the horosphere in our
Cayley-Klein model of $\mathbb{H}^3$:
\begin{equation}
-2 s x^0 x^0-2 x^3 x^3+ 2 (s+1)(x^0 x^3)+(s-1)(x^1 x^1+ x^2 x^2)=0  \tag{3.8}
\end{equation}
\begin{remark}

We have obtained the equation of the horosphere in the Cartesian
coordinate system:
($x:=\frac{x^1}{x^0},~y:=\frac{x^2}{x^0},~z:=\frac{x^3}{x^0}$)
\begin{equation}
\frac{2(x^2+y^2)}{1-s}+\frac{4(z-\frac{s+1}{2})^2}{{(1-s)}^2}=1  \tag{3.9}
\end{equation}
\end{remark}

\begin{remark}

It is useful for visualization purposes to convert the horosphere
equation into polar coordinates. By multiplying the polar coordinate
form by rotation matrices we can easily obtain horospheres around
arbitrary points at infinity in the model. The polar form is by parameters
$\phi\in[0, 2\pi),~\theta\in[0,\pi]$
\begin{equation}
\begin{gathered}
x=\sqrt{\frac{1-s}{2}}\sin{\theta} \cos{\phi}~~~~
y=\sqrt{\frac{1-s}{2}}\sin{\theta} \sin{\phi} \\
z=\frac{1+s}{2}+\frac{1-s}{2}\cos\theta .\tag{3.10}
\end{gathered}
\end{equation}
\end{remark}

\subsection{Volumes of horoball sectors}

The length $l(x)$ of a horospheric arc of a chord segment $x$ is determined
by the classical formula due to {\sc{j.~bolyai}}:
\begin{equation}
l(x)=k \sinh{\frac{x}{k}} \ (\text{at present} \ k=1).  \tag{3.11}
\end{equation}

The intrinsic geometry of the horosphere is Euclidean, therefore,
the area $\mathcal{A}$ of a horospherical triangle is computed by
the formula of Heron. The volume of the horoball pieces can be
calculated using another formula by {\sc{J.~Bolyai}}. If the area of a domain on the horoshere is $\mathcal{A}$, 
the volume determined
by $\mathcal{A}$ and the aggregate of axes drawn from $\mathcal{A}$ is equal to
\begin{equation}
V=\frac{1}{2}k\mathcal{A} \ \ (\text{we assume that} \ k=1 \ \text{in this paper}). \tag{3.12}
\end{equation}

\section{Horoball packings for totally asymptotic Coxeter honeycombs}

\subsection{Basic results on packing in $\mathbb{H}^3$}

In this section we determine optimal horoball packings for the four
totally asymptotic Coxeter tilings $\mathcal{P}_{pqr}$ with
Sch\"alfi symbols $(p,q,r)=(3,3,6),$  $(3,4,4),$  $(4,3,$ $6),$ 
$(5,3,6)$. Vertices of a regular cell $P_{pqr}$ are denoted by $ E_i,
\ (i=0,1,2,3,4 \dots)$, and lie on the absolute of $\mathbb{H}^3$,
hence these vertices are centers of horoballs. The number of the
vertices of $P_{pqr}$ is denoted by $N_{pqr}$ and we write $B_i$ for
the horoball centered  at $E_i$ . We require a horoball to lie at
every ideal vertex of the honeycomb, and we vary the touching types of the horoballs.
\begin{remark}

For example, if $(p,q,r)=(3,3,6)$ then we obtain the ``tetrahedral
case" where $N_{336}=4$. For parameter $(p,q,r)=(3,4,4)$ we get the
``octahedral case"  and $N_{344}=6$.
\end{remark}

The type of a horoball is allowed to expand until either the horoball
comes into contact with other horoballs or a non-adjacent faces of the
honeycomb. These conditions are satisfactory to ensure
that the balls form a non-overlapping horoball arrangement, as such the collection of all
horoballs is a well defined packing in $\mathbb{H}^3$, denoted by
$\mathcal{B}_{pqr}$.

\begin{definition}
The density of a horoball packing in Coxeter honeycomb
$\mathcal{P}_{pqr}$ is defined as
\begin{equation}
\delta(\mathcal{B}_{pqr})=\frac{\sum_{i=1}^{N_{{pqr}}} Vol(B_i \cap P_{pqr})}{Vol(P_{pqr})}. \notag
\end{equation}
\end{definition}

The aim of this section is to determine the optimal packing
densities for the four totally asymptotic tilings (see Table 1) in 3-dimensional hyperbolic space $\mathbb{H}^3$.

We will make heavy use of the following Lemma \cite{Sz05-1}:

\begin{lemma}

Let $\tau_1$ and $\tau_2$ be two congruent trihedra centered at $C_1$ and $C_2$, which share the common edge $C_1C_2$. Let $I(x)$ denote a point along $C_1C_2$ defined as follows; let $B_1(x)$ and $B_2(x)$ be two tangent horoballs centered at $C_1$ and $C_2$ respectively. Pick $I(0)$ to be the point of contact such that the
following equality holds for volumes of horoball sectors:
\begin{equation}
V(0):= 2 Vol(B_1(0) \cap \tau_1) = 2 Vol(B_2(0) \cap \tau_2), \notag
\label{szirmai-lemma}
\end{equation}
and $I(x)$ be the point obtained by the displacement of $I(0)$ by $x$ along $C_1C_2$. 

When $x$ denotes the hyperbolic displacement between $I(0)$ and $I(x)$,
then the volume function
\begin{equation}
V(x):= Vol(B_1(x) \cap \tau_1) + Vol(B_2(x) \cap \tau_2) \notag
\end{equation}
strictly increases as $I(x)$ moves continuously away from $I(0)$.
\label{lemma1}

\end{lemma}
\begin{figure}
\begin{center}
\includegraphics[width=8cm]{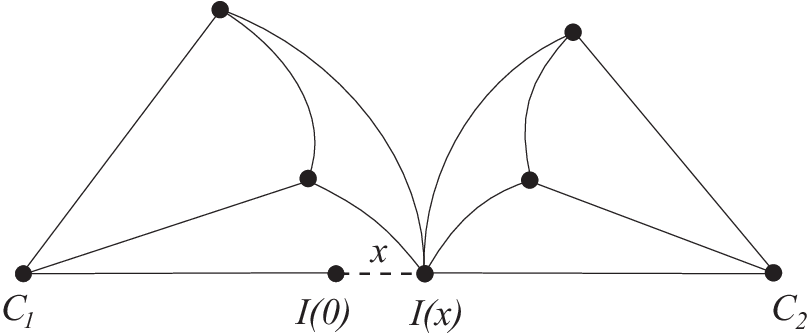}
\caption{$C_1$ and $C_2$ are centers of the horoballs $B_1$ and
$B_2$, and $I(x)$ is the point of contact. $I(0)$ represents the
point at which $Vol(B_1)=Vol(B_2)$.}
\end{center}
\end{figure}
\textbf{Proof:} Let $\mathcal{L}$ and $\mathcal{L}'$ be parallel horocycles with center $C$
and let $A$ and $B$ be two points on the curve $\mathcal{L}$ and $A':=CA \cap \mathcal{L}'$,
$B':=CB \cap \mathcal{L}'$.
By the classical formula of {\sc J. Bolyai} 
$$
\frac{\mathcal{H}(A' B')} {\mathcal{H} (AB)}=e^{\frac{x}{k}}, 
$$
where the horocyclic distance between $A$ and $B$ is denoted by $\mathcal{H} (A,B)$. 

Then by the above formulas we obtain the following volume function:
$$
\aligned
 V(x)& =Vol (B_1 (x) \cap \tau_1)+Vol(B_2(x) \cap \tau_2)= \\
 & =\frac{1}{2} V(0)\left(e^{\frac{2x}{k}}+\frac{1}{e^{\frac{2x}{k}}}\right)=V(0)\cosh\left(\frac{2x}{k}\right).
\endaligned 
$$
It is well known that this function strictly increases in the interval $(0,\infty)$. 
This Lemma is illustrated in Fig.~2. $\square$

\begin{corollary}
Lemma \ref{lemma1} holds true for horoballs intersecting two arbitrary congruent frames of rays joining $C_1$ and $C_2$ as above, respectively.
\end{corollary}

\textbf{Proof:} Follows from the Lemma \ref{lemma1}, by dividing the frames of rays into congruent trihedra. $\square$\\

\subsection{The $(3,3,6)$ tetrahedral tiling}

The $(3,3,6)$ Coxeter tiling is a three
dimensional honeycomb with cells comprised of fully asymptotic
regular tetrahedra. We arbitrarily select one such tetrahedron $E_0E_1E_2E_3$
(see Fig.~3), and place the horoball centers at vertices $E_0, \dots,
E_3$. We vary the types of the horoballs so that they satisfy our constraints of non-overlap.
The packing density is obtained by Definition 4.2.

Define the orthoscheme $A_0A_1A_2A_3$ as follows: $A_0=E_0$ and
$A_3=E_3$ are two vertices of the tetrahedron (see Fig.~3); $A_2$ is
the center of the triangle $E_0E_1E_2$ opposite the vertex $A_3$, and $A_1$ is the
footpoint of $A_3$ on the edge $E_0E_1$. One tetrahedral
cell is decomposable into 6 such congruent orthoschemes. The
Schl\"afli symbol of orthoscheme $A_0A_1A_2A_3$ is $(3,6,3)$, and
the orthoscheme is labeled by $\mathcal{O}_{(3,6,3)}$.

$B_0$ and $B_3$ are two horoballs centered at $E_0$ and $E_3$, i.e.,
the two vertices of the tetrahedron common with the orthoscheme. The
density of the $(3,3,6)$ Coxeter tiling is obtained using 
Definition 4.2:
\begin{equation}
\delta(\mathcal{B}_{336})=\frac{Vol(B_0\cap\mathcal{O}_{(3,6,3)})+Vol(B_3\cap\mathcal{O}_{(3,6,3)})}{Vol({\mathcal{O}_{(3,6,3)}})}
\notag
\end{equation}
\begin{proposition}
The packing density obtained in $\mathcal{O}_{(3,6,3)}$ can be
extended to tetrahedron $P_{336}$ and therefore to the entire
$\mathbb{H}^3$.
\end{proposition}
\textbf{Proof:} We consider the following steps:
\begin{enumerate}
\item In an optimally dense packing, at least two horoballs must touch each other in the tetrahedron, otherwise the density could be improved by blowing up any one horoball until it touches a neighboring horoball.
\item If two horoballs touch at the ``midpoint" of edge $A_0A_3$ as projection of the simplex center on it, then by blowing up the remaining two horoballs,
they will also touch at ``midpoints", due to symmetry considerations.
Note, that this case is the arrangement in the B\"or\"oczky--Florian density upper bound using the same horoballs.
\item Given two horoballs tangent at a non-midpoint of an edge, then the horoball, having the midpoint in its interior, will contain the ``midpoint"
of all $3$ edges extending from its center. As a result the remaining 3 horoballs should be of the same touching type.
The ``small horoballs"  each is tangent only to the large horoball.
\end{enumerate}
We just showed that the ``largest horoball" determines the
configuration of all other horoballs, and as a consequence the
packing density. One parameter corresponding to the ``largest
horoball" suffices to determine the packing density for all
candidates of optimal density.

Due to symmetry, it is enough to consider cases within orthoscheme
$\mathcal{O}_{(3,6,3)}$, where the horoball at $E_0$ expands from the midpoint until it becomes tangent to the side of
the cell opposite to it. Assume no balls cover the midpoint along
$E_0E_3$. Then the horoball $B_0$ can be expanded until the
midpoint. In this case, the horoball $B_3$ is contained within
$\mathcal{O}_{(3,6,3)}$. Finally, the packing density can be varied
by expanding the horoball $B_{0}$ while keeping $B_3$ tangent to it.
If we expand $B_0$ until it touches $A_2$, all candidates for
optimal packings within the orthoscheme are considered. The densities
obtained from the orthoscheme $\mathcal{O}_{(3,6,3)}$ will cover all
candidates for optimally dense packings of honeycomb $(3,3,6)$.
Densities determined within the orthoschemes can be generalized to
the entire $\mathbb{H}^3$ by the symmetries of $\mathcal{P}_{336}$. $\square$

In the rest of this section, we prove the basic theorem on the
optimal packing density in $(3,3,6)$. First, we define the tangent
point $I(0) \in A_0A_3$ of horoballs $B_0$ and $B_3$ so that the
following equality holds for the volumes of the horoball sectors
\begin{equation}
V(0):= 2 Vol(B_0(0) \cap \mathcal{O}_{(3,6,3)}) = 2 Vol(B_3(0) \cap \mathcal{O}_{(3,6,3)}). \notag
\end{equation}

Note that $I(0)$ is not the midpoint of edge $A_0A_3$. Consider point of tangency $I(x)$ of horoballs $B_0(x)$ and $B_3(x)$ along edge $A_0A_3$. Let $x$ denote the hyperbolic displacement of $I(0)$ to $I(x)$. The volume function
$V(x)$ is defined as follows:
\begin{equation}
V(x):= Vol(B_0(x) \cap \mathcal{O}_{(3,6,3)}) + Vol(B_3(x) \cap \mathcal{O}_{(3,6,3)}) .\notag
\end{equation}

By Lemma \ref{lemma1} it follows that function $V(x)$ strictly increases as
$I(x)$ ($x \in [-\text{arctanh}(1/2), \text{arctanh}(1/2)]$) moves away from $I(0)$ along $A_0A_3$.
That implies that the density function %

\begin{equation}
\delta(\mathcal{B}_{336})(x)=\frac{V(x)}{\mathcal{O}_{(3,6,3)}}=\frac{Vol(B_0(0) \cap \mathcal{O}_{(3,6,3)})e^{2x} + Vol(B_3(0) \cap \mathcal{O}_{(3,6,3)})e^{-2x}}{Vol(\mathcal{O}_{(3,6,3)})} \notag
\end{equation}
attains its maximum at the two endpoints of the interval $[-\text{arctanh}(1/2), \text{arctanh}(1/2)]$.

\begin{definition}
A \textit{coordinate system} is assigned to the orthoscheme
$A_0A_1A_2A_3$; let $A_2:=(1,0,0,0)$ be the origin,
$A_0:=(1,0,1,0)$, $A_3:=(1,0,0,1)$ and
$A_1=(1,\frac{\sqrt{3}}{4},\frac{1}{4},0)$ where $A_1$ is the
``midpoint" of the edge $E_0E_1$ (see Fig.~3).
\end{definition}
\begin{figure}
\begin{center}
\includegraphics[width=11cm]{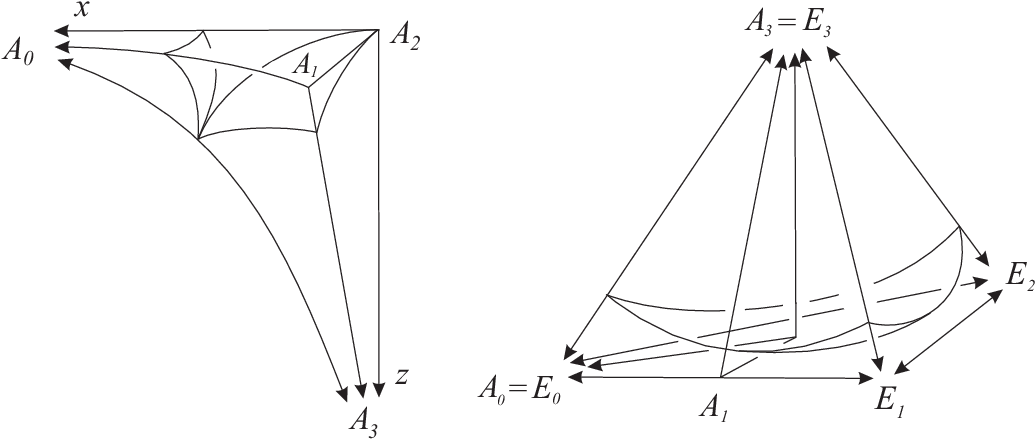}
\caption{a. Orthoscheme $\mathcal{O}_{(3,6,3)}$; b. the totally asymptotic regular tetrahedron $E_0E_1E_2E_3$ for $(3,3,6)$.}
\end{center}
\end{figure}

In the orthoscheme $A_0A_1A_2A_3$, the horoball centers are located
at points $A_0$ and $A_3$, and the horoballs meet along edge
$A_0A_3$ (See Fig.~3 (i)). In order to determine the properties of
the packing, we calculate the five points of intersection of the
horospheres with the edges of the orthoscheme. The distances of these
five points of intersection determine the area of the horospheric
triangles, hence the volume of the horoball sectors through Bolyai's
formulas (3.11-12). 

To aid the discussion of our computations regarding horospheres we introduce the $s$-parameter as the following. 
Define $s$ to be the Euclidean distance between the origin $A_2$ and point $S(1,0,0,s)$ along the axis defined by $A_2A_3$.
Based on the proof of Proposition 4.4, $s$ is the only parameter that determines the
density of the packing.

All horoball configurations that yield valid candidates for the
optimal packing density occur while we continuously vary the
horoball s-parameter on $[0, 1/2]$. In order to determine the
optimal packing, we proceed as follows. We express the density
$\delta(\mathcal{B}_{336})$ as a function of $s$, study its
behavior and determine its extremal points. If $s=1/2$, all
horoballs are in same type, thus we achieve the packing arrangement
$\mathcal{B}_{336}^1$, which is the B\"or\"oczky--Florian case with
the known packing density. On the other hand, for $s=0$, we find a
different horoball arrangement $\mathcal{B}_{336}^2$ (see Fig.~3)
with the same density as for $s=1/2$. Finally, from Lemma \ref{lemma1} it follows that all
other densities with horoball parameter $s \in (0,1/2)$ are smaller.

First we calculate the five intersections of the edges of
$\mathcal{O}_{(3,6,3)}$ and the two horoballs $B_0$ and $B_3$ for
$s=0$. Recall that these are all a function of the  ``type" of the
large horoball, and depend on the parameter $s$. The length of the
sides of the horospherical triangle, with vertices $X_i,Y_i,Z_i \in
B_i$, $i\in\{0,3\}$ are given by
\begin{equation}
a_i:=2\sinh\frac{d(X_i,Y_i)}{2}~~~~~
b_i:=2\sinh\frac{d(Y_i,Z_i)}{2}~~~~~
c_i:=2\sinh\frac{d(Z_i,X_i)}{2}. \notag
\end{equation}

Based on the side lengths, Heron's formula can be used to obtain the
area $\mathcal{A}_i(s)$ of horospherical triangles
$X_iY_iZ_i$. Using Bolyai's volume formula (3.12) for
horoball pieces, the Definition 4.2, and setting $k=1$, the
density of the packing can be expressed as a function of $s$ as well
\begin{equation}
\delta(\mathcal{B}_{336})(s)=\frac{\frac{1}{2}(\mathcal{A}_0(s)+\mathcal{A}_3(s))}{Vol({\mathcal{O}_{(3,6,3)}})}, \ \ s \in [0,1/2]. \notag
\end{equation}

From Lemma \ref{lemma1} it follows, that the optimal densities are realized on the endpoints of the interval $[0,1/2]$.
In order to determine the highest packing density, we calculate the
density for horoball pieces with $s=0$. The volume of orthoscheme
$\mathcal{O}_{(3,6,3)}$ in $\mathbb{H}^3$ is calculated using
Lobachevsky's volume formula (see \cite{Sz03-1}). All other
densities with $ 0 \leq s \leq 1/2$ can be evaluated using the
corresponding volume formula. The results are displayed in Fig.~4.a,
while Fig.~4.b and Fig.~4.c illustrate the ball arrangements
$\mathcal{B}_{336}^i \ \ (i=1,2)$ for $s=1/2$ and $s=0$,
respectively. Note that $\mathcal{B}_{336}^1$ and
$\mathcal{B}_{336}^2$ are two different cases. In the first case the
horoballs are of the same type and they touch each other pairwise, while
for $s=0$ we have balls of two distinct types, and the ``small balls"
do not meet. As a result, we have just proved: %
\begin{figure}
\begin{center}
\includegraphics[width=4cm]{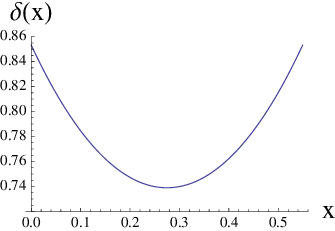}
\includegraphics[width=4cm]{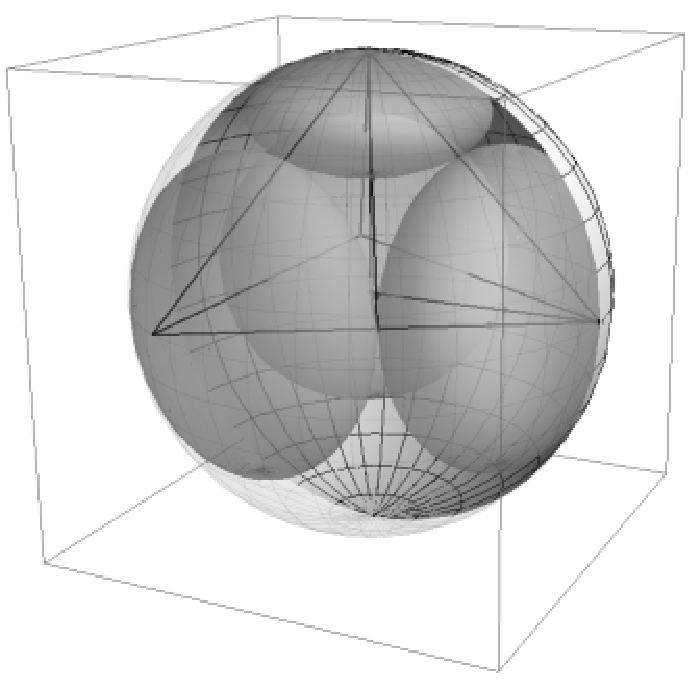}
\includegraphics[width=4cm]{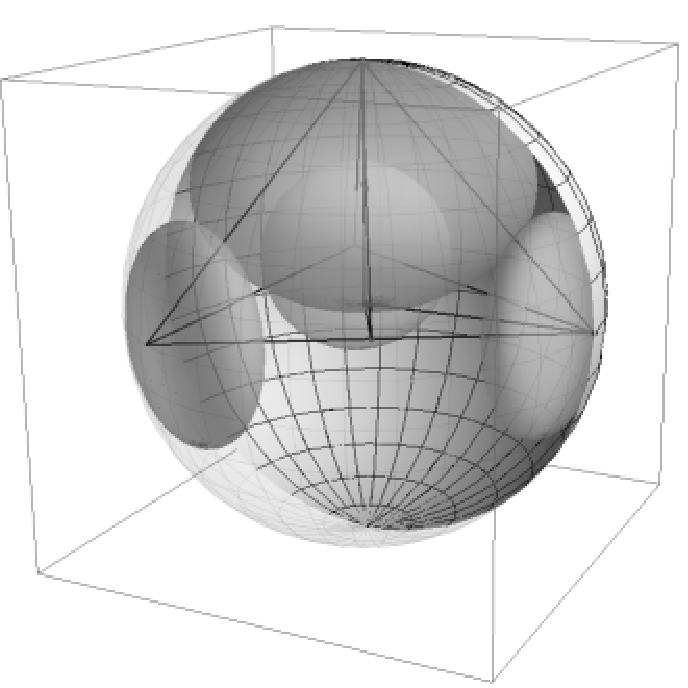}
\\a. $s\in (0,1/2)$~~~~~~~~~~~~~~~~~~~~~b. $s=1/2$ ~~~~~~~~~~~~~~~~~~~~~~c. $s=0$ \\
\caption{Horoball packings of $(3,3,6)$ tiling.}
\end{center}
\end{figure}
\begin{theorem}
There are two distinct optimally dense horoball arrangements
$\mathcal{B}_{336}^i,$ $(i=1,2)$ for the tetrahedral Coxeter tiling
$(3,3,6)$ with the same density: $\delta(\mathcal{B}_{336}^i) \approx
0.85327609$.
\end{theorem}

\subsection{The $(3,4,4)$ octahedral tiling}

In this section we consider horoball packings with centers located at ideal vertices of the
octahedral honeycomb $(3,4,4)$. Our approach is as in
the previous section for the $(3,3,6)$ Coxeter tiling. We again
allow the horoballs to be of different types. The $(3,4,4)$ Coxeter
tiling decomposes the 3-dimensional hyperbolic space into congruent
cells consisting of regular fully asymptotic  octahedra. Four
octahedra meet along each edge. As in the case of $(3, 3, 6)$, we again 
choose one cell of the tiling in order to perform our density
calculations. Again, we vary the types of the horoballs centered in
the honeycomb vertices. The packing density obtained form this
cell is again extended to $\mathbb{H}^3$ by $Sym{\mathcal{P}}_{344}$.

We define the orthoscheme $A_0A_1A_2A_3$ of an octahedron for the
calculations. 

\begin{definition}
The following \textit{coordinate system} is defined for orthoscheme\\
$A_0A_1A_2A_3$ (see Fig. 5):
$$A_2=(1,0,0,0), A_0=(1,0,1,0), A_3=(1,0,0,1), A_1=(1,\frac{1}{2},\frac{1}{2},0).$$
\end{definition}

Let $A_0$ and $A_3$ be two adjacent vertices of the
octahedron $E_0E_1E_2E_3E_4E_5$ (see Fig.~5.~a-b-c), $A_2$ be the
center of the octahedron, and take $A_1$ as the midpoint, in the
Euclidean sense, of the edge extending from $A_0$ sharing a common
facet with $A_3$. This orthoscheme is $\mathcal{O}_{444}$ and has
Schl\"afli symbol $(4,4,4)$.

Using the Lobachevsky volume formula for orthoschemes, we
obtain the volume of one octahedron \cite{Sz05-2}: $Vol({P}_{344})=
16 \cdot Vol(\mathcal{O}_{(4,4,4)}) \approx 3.66384$. Applying the
definition of the packing density for the case of tiling $(3,4,4)$,
we obtain:

\begin{equation}
\delta(\mathcal{B}_{344})=\frac{\sum_{i=1}^{6} Vol(B_i \cap P_{344})}{Vol(P_{344})}, \tag{4.1}
\end{equation}
where $B_i\cap {P}_{344}$ $(i=1,\dots, 6)$ denote the 6 horoball sectors, one in each vertex of the octahedron $P_{344}$,
and we assume that the horoballs $B_i$ form a horoball packing in $\mathbb{H}^3$.

We consider the following three basic horoball configurations
$\mathcal{B}_{344}^i$, $(i=1,2,3)$:

\begin{enumerate}
\item All 6 horoballs are of the same type and the adjacent horoballs
touch each other at the ``midpoints" of each edge.
We define the point of tangency of two horoballs $B_0$ and $B_3$ on
side $A_0A_3$ to be $I(0)$ so that the following equality holds:
\begin{equation}
V(0):= 6 \cdot  Vol(B_0(0) \cap {P}_{344)}) = 6 \cdot Vol(B_3(0) \cap {P}_{344})=6 \cdot V_0. \notag
\end{equation}
In this case $V_0:=Vol(B_i\cap {P}_{344}) = 0.5$
$(i=1,\dots, 6)$ (see Fig.~5.a, Fig.~6.c).
\item Two ``larger horoballs" with centers at $E_3$ and $E_5$ are
tangent at the center of the octahedron, while horoballs at the remaining six vertices touch both ``larger"
horoballs. The point of tangency of the above two horoball types on
segment $I(0)E_0$ is denoted by $I(x_1)$ where $x_1=\log(2)/2$ is the
hyperbolic distance between $I(0)$ and $I(x_1)$. In this case
$V_1:=Vol(B_i\cap {P}_{344}) = 1$ $(i=3,5)$ and
$V_2:=Vol(B_i\cap {P}_{344}) = 0.25$ $(i=0,1,2,4)$ (see
Fig.~5.b, Fig.~6.b).
\item One horoball of the ``maximally large" type centered at $E_3$.
The large horoball is tangent to all non-neighboring faces of the
octahedron and it determines the other five horoballs touching the
``large horoball". The point of tangency of the two horoballs along
segment $I(x_1)E_0$ is denoted by $I(x_2)$ where $x_2=-\log(2)/2$ is the
hyperbolic distance between $I(0)$ and $I(x_2)$. In this case
$V_3:=Vol(B_3\cap {P}_{344}) = 0.25$, $V_5:=Vol(B_5\cap {P}_{344}) = 0.0625$ and $V_i:=Vol(B_i\cap{P}_{344}) = 0.03125$ $(i=0,1,2,4)$ 
(see Fig.~5.c, Fig.~6.a).
\end{enumerate}
\begin{figure}
\begin{center}
\includegraphics[width=13cm]{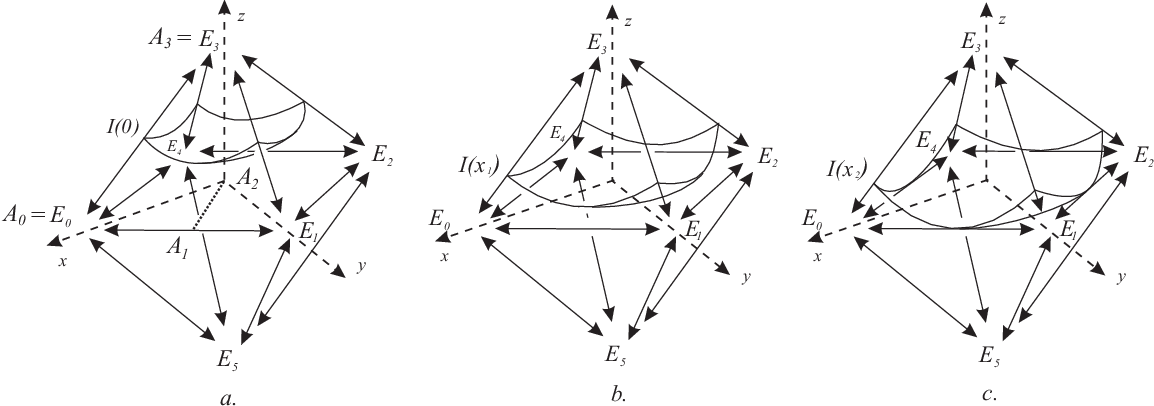}
\caption{The site of the horoball $B_3$ of the basic horoball
arrangements $\mathcal{B}_{344}^i$ $(i=1,2,3)$ in the octahedron
$P_{344}$.}
\end{center}
\end{figure}
Due to symmetry considerations it is sufficient to restrict our attention to the family of
cases in which the type of one horoball is varied between the following two limiting cases: passing through the ``midpoint" of an edge and touching the opposite facet of the cell. Here we give an
analogous argument to that given in the proof of Proposition 4.4. Assume that
none of the horoballs covers a ``midpoint". Then the packing density
may be improved until at least one ball reaches a midpoint.
Moreover, consider that point $I(x)$ is on edge $A_0A_3=E_0E_3$.
This point is where the horoballs $B_i(x), \ (i=0,3)$ are tangent at
point $I(x) \in A_0I(0)$. Then $x$ is the hyperbolic distance
between $I(0)$ and $I(x)$, and it is analogous to the previous section. It is easy to see that we have to study two
different cases to determine the optimal horoball arrangement:
\begin{enumerate}
\item $x \in [0,x_1]$, horoballs $B_3$ and $B_5$ touch horoballs $B_i$ $(i=1,2,3,4)$.
\item $x \in [x_1,x_2]$, horoball $B_3$ touches horoballs $B_i$ $(i=1,2,3,4,5)$.
\end{enumerate}
In the {\bf first case} the function $V(x)$ can be computed by the
following formula
\begin{equation}
V(x):= 4 \cdot Vol(B_0(x) \cap {P}_{344}) + 2 \cdot Vol(B_3(x) \cap {P}_{344}) \ \ x\in [0,x_1]. \notag
\end{equation}
Similarly to the Lemma \ref{lemma1}, we can prove the following Lemma:
\begin{lemma}
\begin{equation}
\begin{gathered}
V(x):= 4 Vol(B_0(x) \cap {P}_{344}) + 2 Vol(B_3(x) \cap {P}_{344})= \\
=V_0 \left(2 e^{2 x} + 4 e^{-2x}\right), \ \ \ x \in [0,x_1], \notag
\end{gathered}
\end{equation}
and the maxima of function $V(x)$ are realized in points $I(0)$ and $I(x_1)$.
\label{lemma2}
\end{lemma}

In the {\bf second case}, similarly to Lemma \ref{lemma2}, the volume
function $V(x)$ is given by the following formula: %
\begin{lemma}
\begin{equation}
\begin{gathered}
V(x):= V_1 \left(e^{2 (x-x_1)} +  e^{-2(x-x_1)}\right) + V_2 \left( 4 e^{-2(x-x_1)}\right), \ \ \ x \in [x_1,x_2], \notag
\end{gathered}
\end{equation}
and the maxima of function $V(x)$ is achieved at points $I(x_1)$
and $I(x_2)$.
\label{lemma3}
\end{lemma}

By Definition 4.2 and Remark 4.1, as well as by Lemmas \ref{lemma2}-\ref{lemma3}, it
follows that the densest horoball packing is realized in three
distinct primary  horoball arrangements $\mathcal{B}_{344}^i$
$(i=1,2,3)$. These optimal horoball packings belong to horoball
parameters $s=-1/3$, $s=0$, $s=1/3$ (see Section (3.1)) and are
illustrated in Fig.~6.a-b-c.

The maximal density can be computed by the method described in
Section 4.2. Thus we have proven %
\begin{theorem}
Three different optimally dense horoball arrangements
$\mathcal{B}_{344}^i,$ $(i=1,2,3)$ exist for the octahedral Coxeter
tiling $(3,4,4)$, which share the optimum density
$\delta(\mathcal{B}_{344}^i) \approx 0.818808$.
\end{theorem}
\begin{figure}[h]
\begin{center}
\includegraphics[width=4cm]{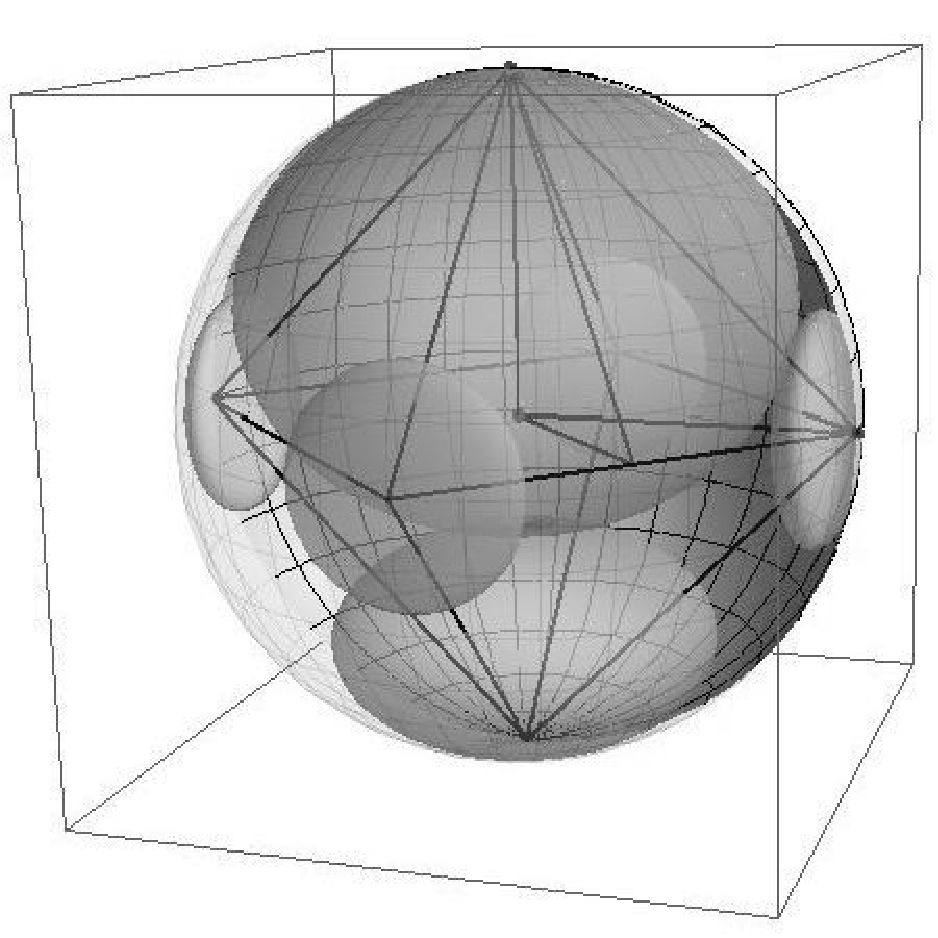}
\includegraphics[width=4cm]{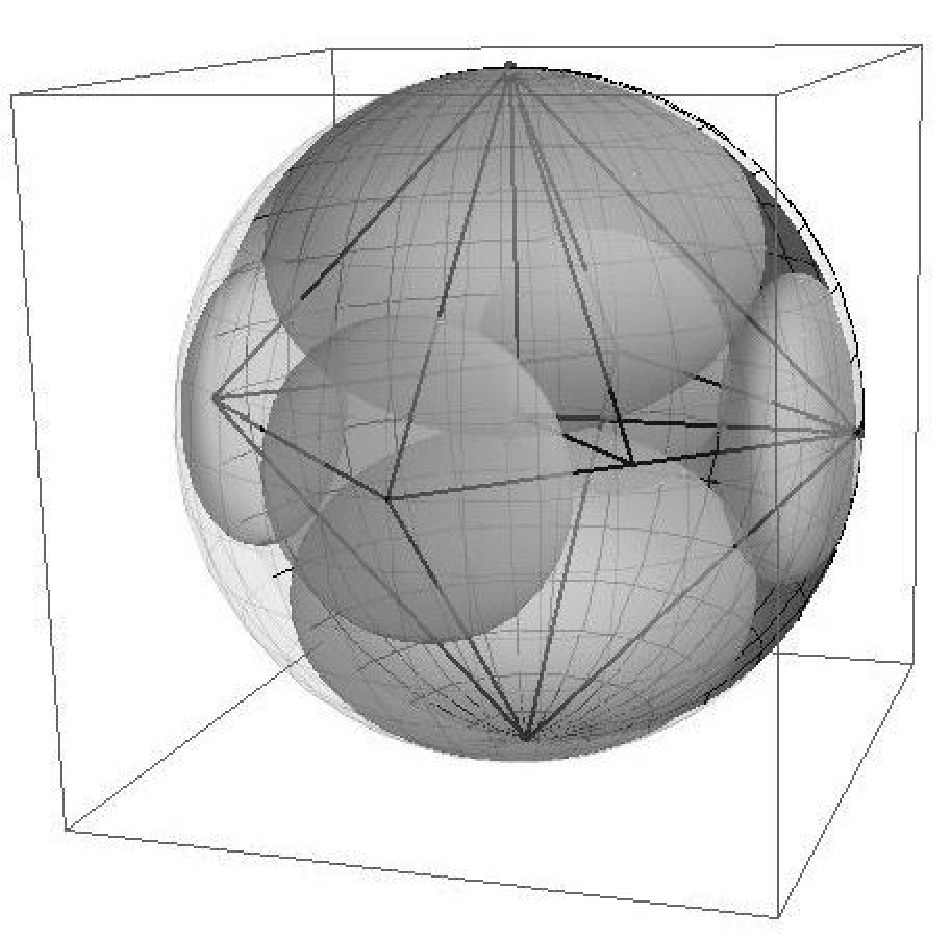}
\includegraphics[width=4cm]{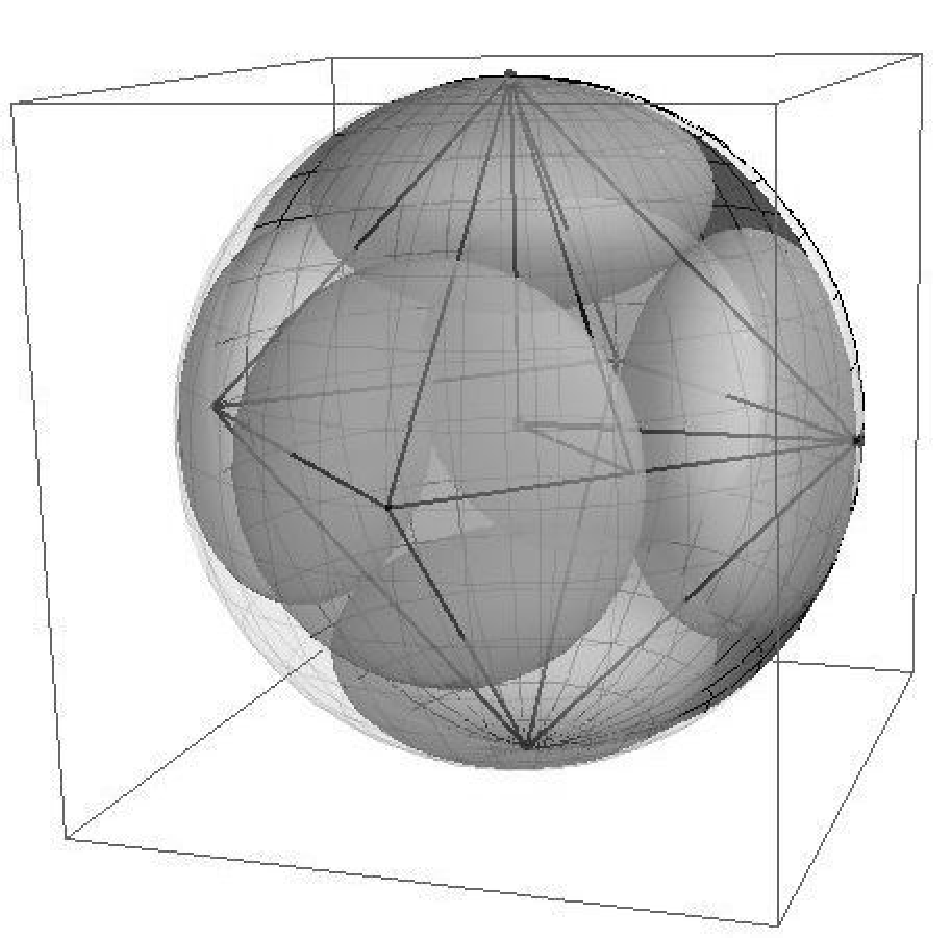}
\\a. $s=-1/3$ ~~~~~~~~~~~~~~~~~~~~~b. $s=0$~~~~~~~~~~~~~~~~~~~~~~c. $s=1/3$ \\
\caption{Optimal horoball packings of Coxeter honeycomb $(3,4,4)$.}
\end{center}
\end{figure}

\subsection{The $(4,3,6)$ Cubic Tiling}

The optimal packing densities for the cubic Coxeter tiling $(4,3,6)$
can be obtained by similar approaches as in the previous two
sections for $(3,3,6)$ and $(3,4,4)$. 

We consider horoball packings with centers located at the ideal vertices of the
cube honeycomb $(4,3,6)$. 

Analogous to the above cases we introduce a projective coordinate system, 
by an orthogonal vector basis with signature $(-1,1,1,1)$, 
with the following coordinates of the vertices of the infinite 
regular cube (see Fig.~7), in the Cayley-Klein ball model: 
$$
E_0(1,-\frac{\sqrt{2}}{\sqrt{3}},\frac{\sqrt{2}}{3},\frac{1}{3}), \ 
E_1(1,-\frac{\sqrt{2}}{\sqrt{3}},-\frac{\sqrt{2}}{3},-\frac{1}{3}), \
E_2(1,0,2\frac{\sqrt{2}}{3},-\frac{1}{3}),
$$
$$
E_3(1,0,0,1), \
E_4(1,\frac{\sqrt{2}}{\sqrt{3}},-\frac{\sqrt{2}}{3},-\frac{1}{3}).
$$
Using the Lobachevsky volume formula for orthoschemes, we
obtain the volume of one cubic cell \cite{Sz05-2}: $Vol({P}_{436})=
48 \cdot Vol(\mathcal{O}_{(4,3,6)}) \approx 5.07471$. Applying the
definition of the packing density for the case of tiling $(4,3,6)$,
we obtain:
\begin{equation}
\delta(\mathcal{B}_{436})=\frac{\sum_{i=1}^{8} Vol(B_i \cap P_{436})}{Vol(P_{436})}, \tag{4.7}
\end{equation}
where $B_i\cap {P}_{436}$ $(i=1,\dots, 8)$ denote the 8 horoball sectors, one in each vertex of the cube $P_{436}$ and we assume that the 
horoballs $B_i$ form a packing in $\mathbb{H}^3$.
\begin{figure}
\begin{center}
\includegraphics[width=13cm]{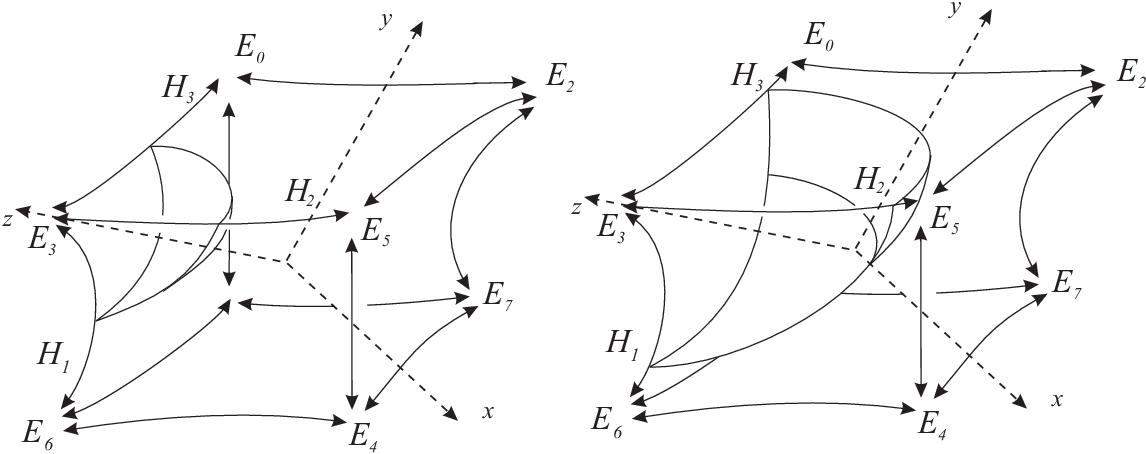}
\caption{The site of the horoball $B_3$ of the base horoball
arrangements $\mathcal{B}_{436}^i$ $(i=1,2,3,\dots,8)$ in the cubic tiling
$P_{436}$.}
\end{center}
\end{figure}
We consider the following four main horoball configurations
$\mathcal{B}_{436}^i$ $(i=1,2,3,4)$:

\begin{enumerate}
\item All 8 horoballs are of the same type and the adjacent horoballs
touch each other at the ``midpoints" of each edge (see Fig.~8.a).
The density of this packing:
$\delta(\mathcal{B}_{436}^1)\approx 0.682621$.
\item Two ``larger horoballs" with centers at $E_3$ and $E_7$ are
tangent at the center of the cube, while the congruent horoballs at the
the remaining six vertices touch one of their respective neighboring ``larger
horoballs" (see Fig.~8.b). 
The density of this packing: $\delta(\mathcal{B}_{436}^2)\approx 0.682621$.
\item The horoballs centered at the two complementary tetrahedral sublattices of the cube are of the same type respectively. Horoballs centered on one sublattice touch at the ``midpoints" of the ``diagonals" of cube's facets.
The other four horoballs corresponding to the other sublattice touch the adjacent ``large horoballs" (see Fig.~9.a).
The density of this packing: $\delta(\mathcal{B}_{436}^3)\approx 0.853276$.
\begin{figure}[ht]
\begin{center}
\includegraphics[width=4cm]{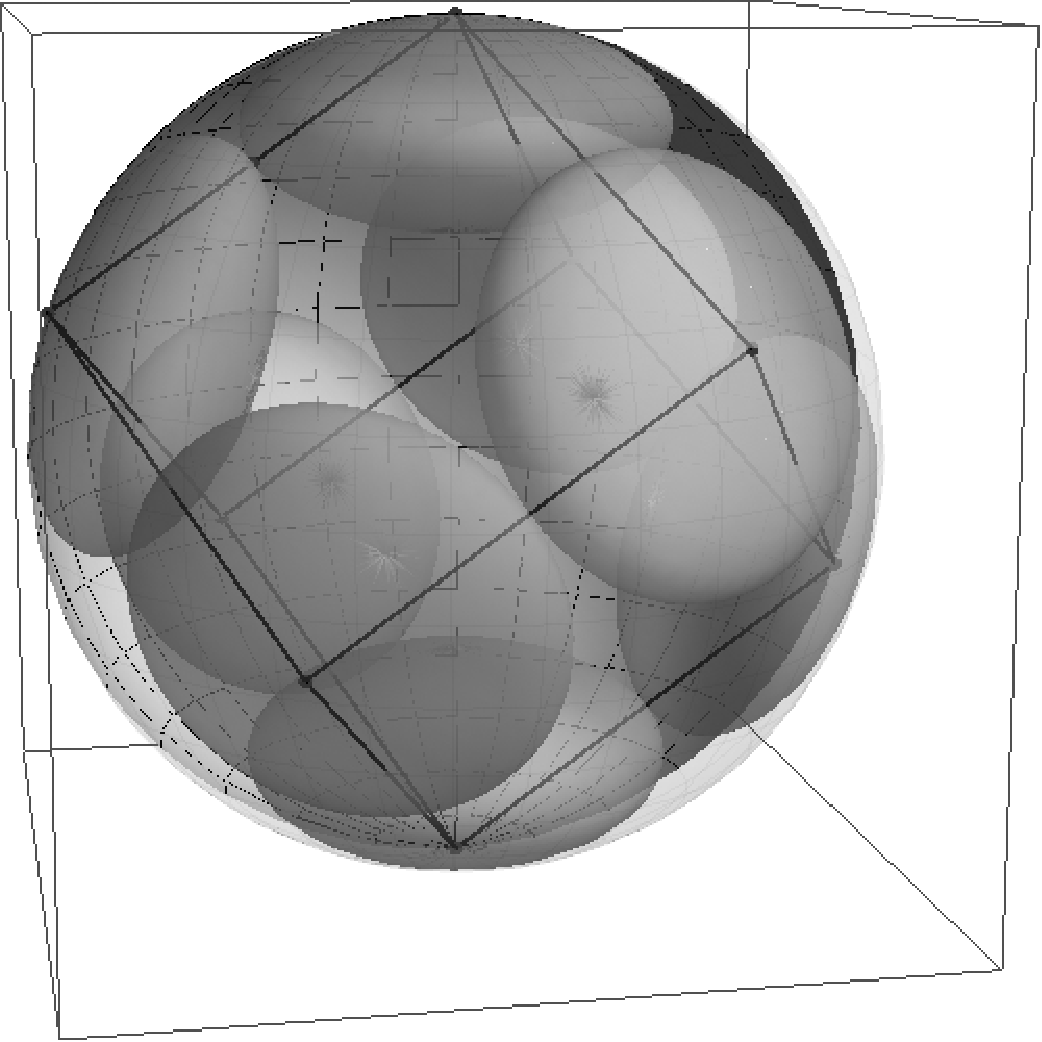} ~~~~~~~~~~~~
\includegraphics[width=4cm]{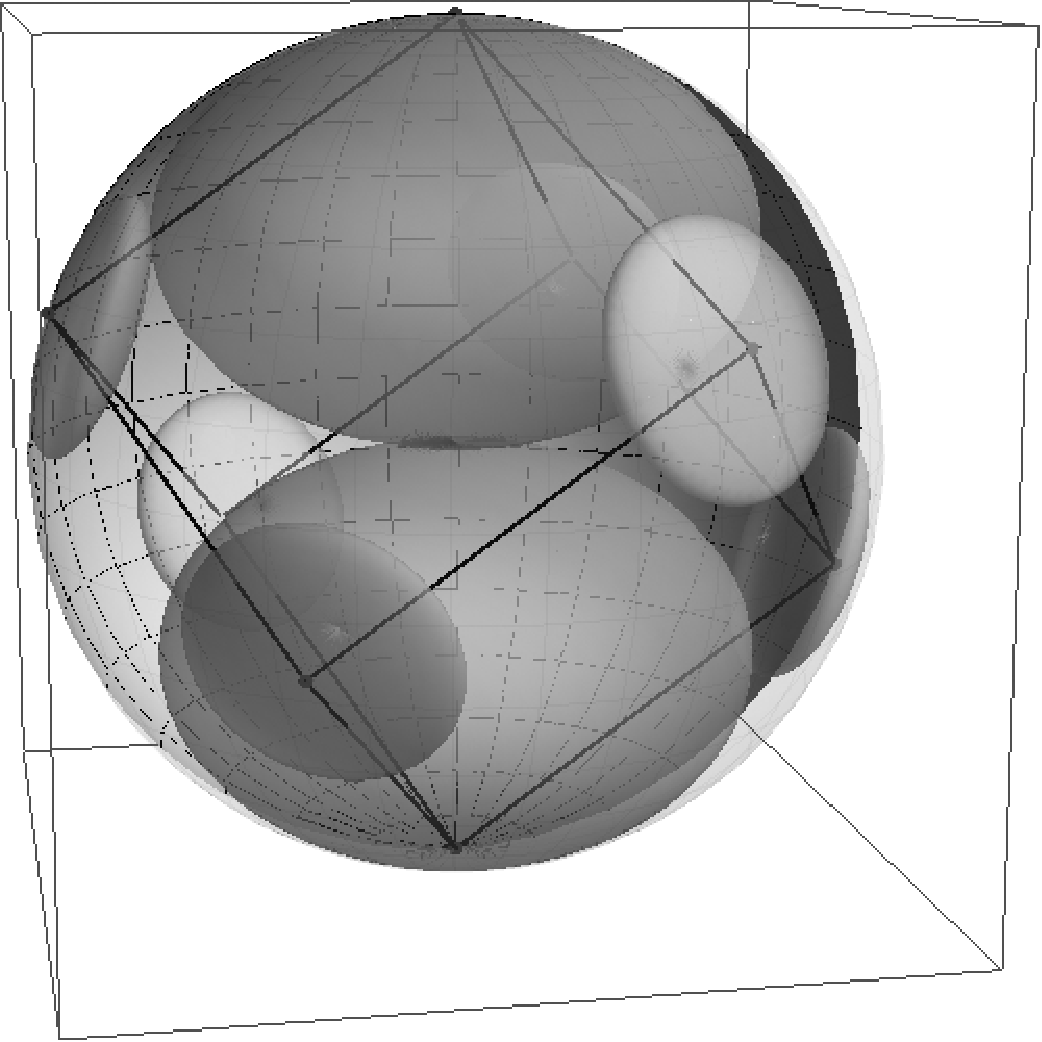}
\\a. ~~~~~~~~~~~~~~~~~~~~~~~~~~~~~~~~~~~~~~~ b. \\
\caption{Locally optimal packings of Coxeter honeycomb $(4,3,6)$}
\end{center}
\end{figure}
\item One horoball of the ``maximally large" type centered at $E_3$.
The large horoball is tangent to all non-neighboring sides of the
cube and it determines the other five horoballs touching the
``large horoball" (see Fig.~9.b).
The density of this packing: $\delta(\mathcal{B}_{436}^4)\approx 0.853276$.
\end{enumerate}
\begin{figure}
\begin{center}
\includegraphics[width=4cm]{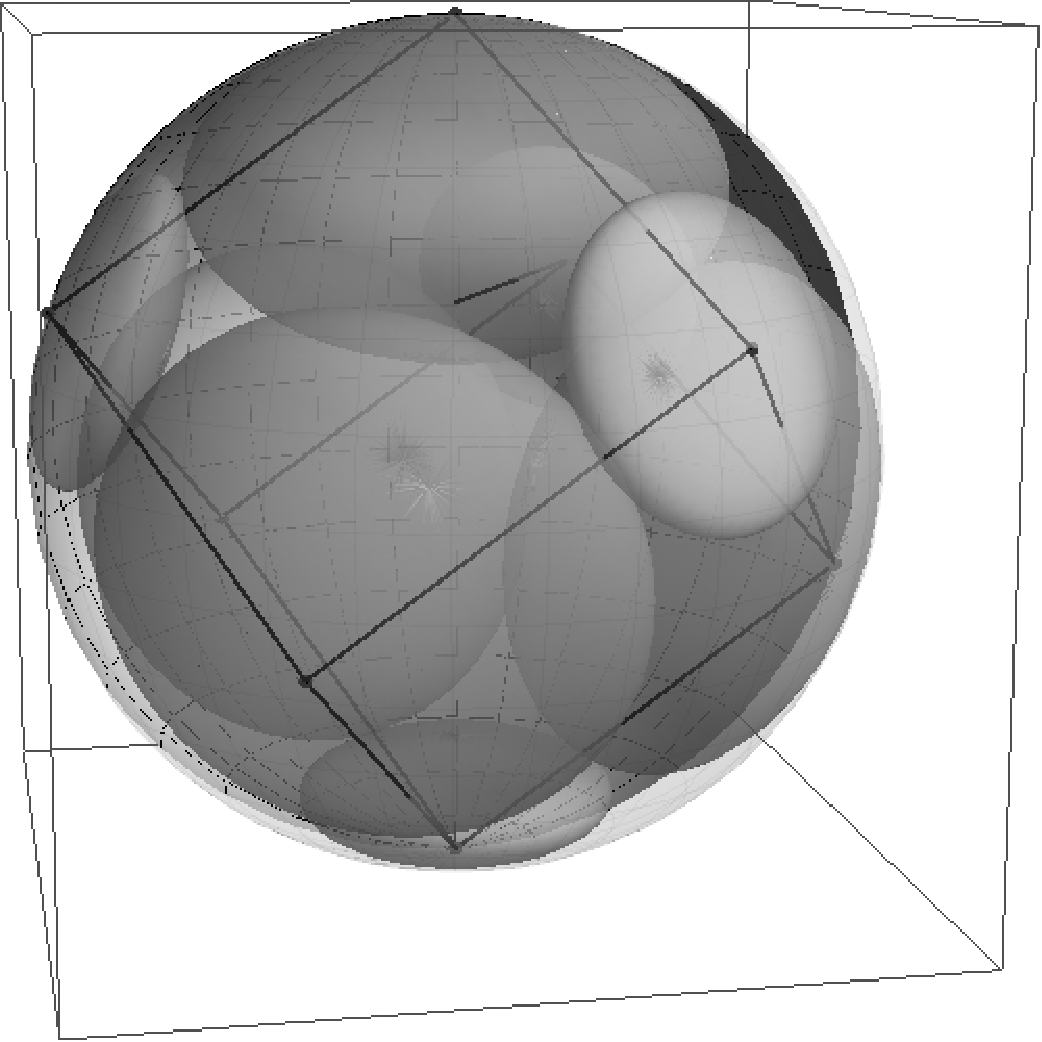} ~~~~~~~~~~~~
\includegraphics[width=4cm]{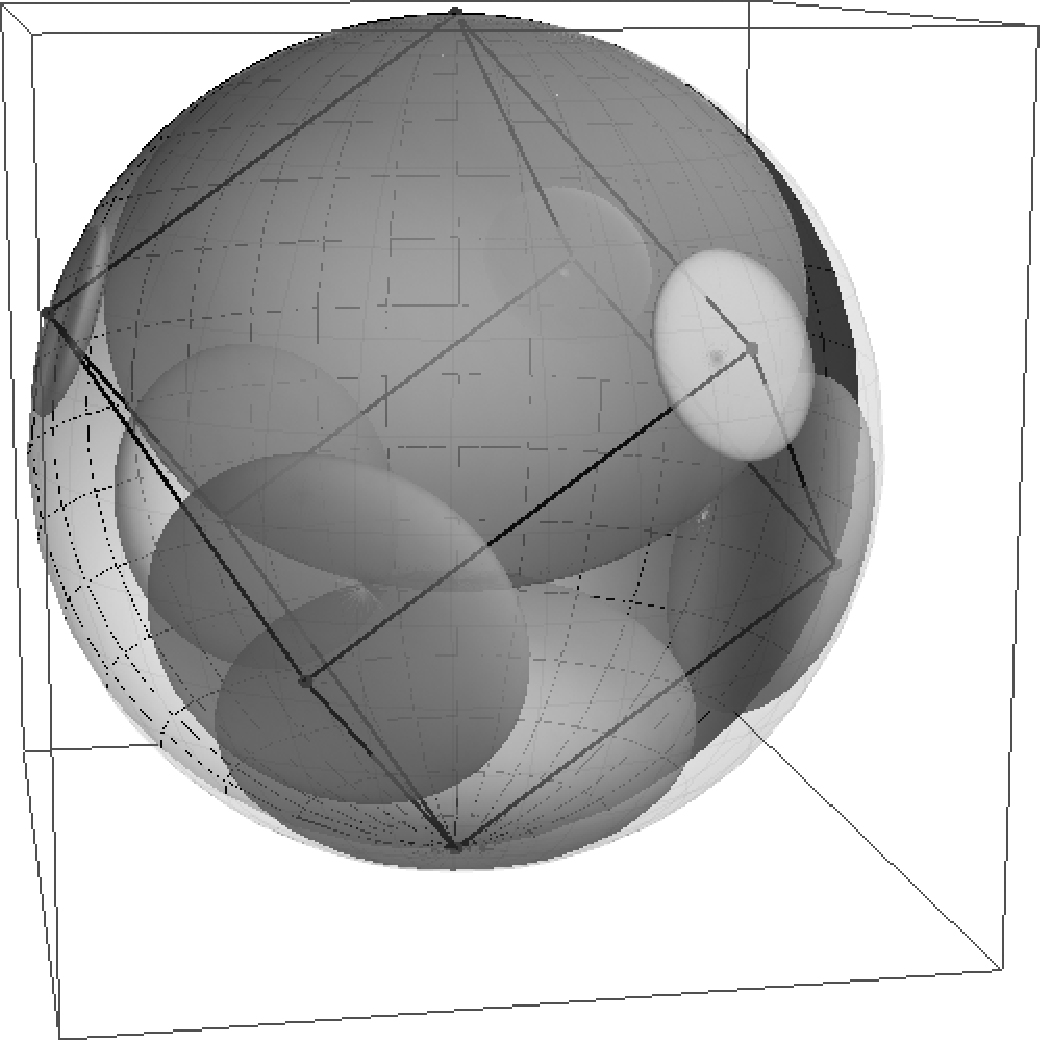}
\\a.  ~~~~~~~~~~~~~~~~~~~~~~~~~~~~~~~~~~~~~~ b.  \\
\caption{Optimal horoball packings of Coxeter honeycomb $(4,3,6)$.}
\end{center}
\end{figure}
By using Lemma \ref{lemma2} and its
corollaries, similarly to the above cases we may again write the volume function of the horoball
pieces to prove the optimality of two of the four limiting cases.  
Finally, we obtain the following
\begin{theorem}
Two different optimally dense horoball arrangements
$\mathcal{B}_{436}^i,$ $(i=3,4)$ exist for the cubic Coxeter
tiling $(4,3,6)$, which share the optimum density
$\delta(\mathcal{B}_{436}^i) \approx 0.85327609$.
\end{theorem}
\section{The $(5,3,6)$ Dodecahedral Tiling}

The optimal packing density for the dodecahedral Coxeter tiling $(5,3,6)$
is obtained through a similar method as in the previous three subsections, hence we omit the details. 

We consider horoball packings with centers of
horoballs located at ideal vertices of the
dodecahedral honeycomb $(5,3,6)$.

Analogous to the previous cases we introduce the projective coordinate system, 
by an orthogonal vector basis with signature $(-1,1,1,1)$, 
with the following coordinates of the vertices of the infinite 
regular dodecahedron, in the Cayley-Klein ball model. The dodecahedron contains a cubic sub-lattice with coordinates adopted from the previous section, as well as 12 other vertices obtainable through rotation about these vertices. 

Using the Lobachevsky volume formula for orthoschemes, we
obtain the volume of one dodecahedron \cite{Sz05-2}: $Vol({P}_{536})=
120 \cdot Vol(\mathcal{O}_{(5,3,6)}) \approx 20.580199\dots$. Applying the
definition of the packing density for the case of tiling $(5,3,6)$,
we obtain:
\begin{equation}
\delta(\mathcal{B}_{536})=\frac{\sum_{i=1}^{20} Vol(B_i \cap P_{536})}{Vol(P_{536})}, \tag{4.7}
\end{equation}
where $B_i\cap {P}_{536}$ $(i=1,\dots, 20)$ denote the 20 horoball sectors, one in each vertex of the dodecahedron $P_{536}$.
By using that the dodecahedral tiling has a cubic sublattice, we consider the following five main horoball configurations
$\mathcal{B}_{536}^i$ $(i=1,\dots, 5)$:

\begin{enumerate}
\item All 20 horoballs are of the same type and adjacent horoballs
are tangent at the ``midpoints" of each connecting edge.
The density of this packing:
$\delta(\mathcal{B}_{536}^1)\approx 0.550841$.
\item Two types of horoballs occur in this packing confguration. Eight larger horoballs are centered at the lattice points of the dodecahedron making 
up a cubic sublattice as in packing $\mathcal{B}_{436}^1$. Twelve congruent horoballs are located at the remaining 12 lattice points. 
The density of this packing: $\delta(\mathcal{B}_{536}^2)\approx 0.70309$.
\item This packing configuration contains horoballs of 4 types. The cubic sublattice within the dodecahedral lattice 
has the same ball configuration as $\mathcal{B}_{436}^2$, and the two types of the balls on the cubic lattice points uniquely determine two types of 
neighboring horoballs.
The density of this packing: $\delta(\mathcal{B}_{536}^3)\approx 0.78725$.
\item There are three types of horoballs in this packing. The cubic sublattice within the dodecahedral lattice has the same ball configuration 
as $\mathcal{B}_{436}^3$, and the larger horoball uniquely determines the type of the remaining 12 horoballs. 
The density of this packing: $\delta(\mathcal{B}_{536}^4)\approx 0.784181$.
\item This limiting case is an extension of packing  $\mathcal{B}_{536}^3$. 
We inflate the horoball located at $(1,0,0,1)$ until it touches the non-adjacent side. 
This packing consists of 6 horoball types, 4 of which are on the cubic sublattice and are a non-limiting case of the cubic tiling. 
These 4 horoballs uniquely determine two horoball types of the remaining 12 vertices of the dodecahedral lattice. 
The density of this packing: $\delta(\mathcal{B}_{536}^5)\approx 0.71246$. 
\end{enumerate}

\begin{figure}[ht]
\begin{center}
\includegraphics[width=6cm]{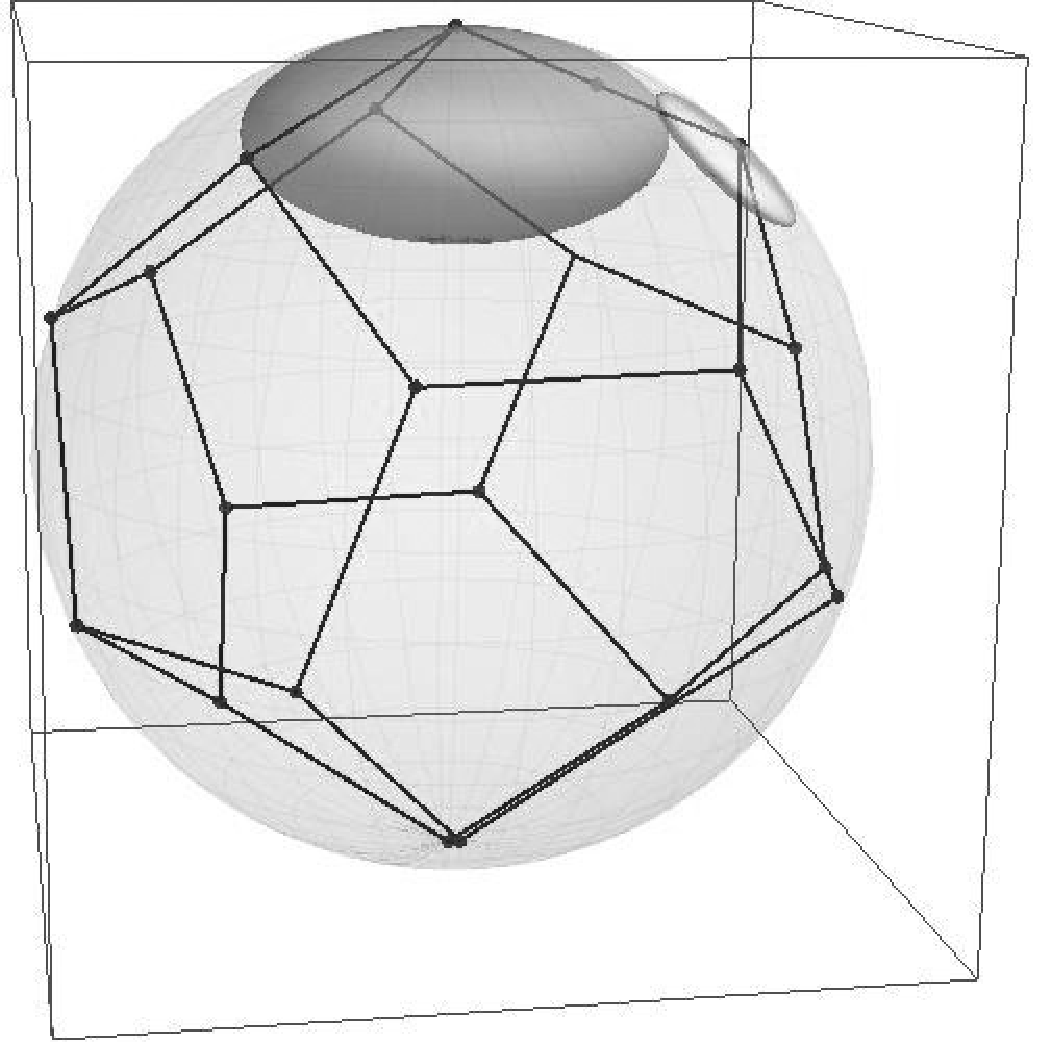} 
\caption{Two tangent horoballs in a packing of the $(5,3,6)$ honeycomb.}
\end{center}
\end{figure}

By using Lemma \ref{lemma2} and its
corollaries, similarly to the above three cases we may again write the volume function of the horoball
pieces to prove the optimality of two of the four limiting cases, leading to the following theorem:

\begin{theorem}
The optimally dense horoball arrangement
$\mathcal{B}_{536}^3$ for the dodecahedral Coxeter
tiling $(5,3,6)$ has optimal density
$\delta(\mathcal{B}_{536}^3) \approx 0.787251$.
\end{theorem}

\section{Conclusion}

Locally optimal horoball packings were studied in this paper,
in which different types of horoballs were placed at the lattice
points of fully asymptotic Coxeter honeycombs. We proved that
the value obtained by B\"or\"oczky and Florian as the universal
density upper bound for all congruent ball packings in hyperbolic $3$-space
remains the upper bound even if horoballs of different types are considered.
Moreover, there are two distinct optimally dense packings for
the $(3,3,6)$ Coxeter tiling. We again encounter the B\"or\"oczky -- Florian density
upper bound as upper limits when varying the types of horoballs at
the lattice points of the $(4,3,6)$ cubic Coxeter tilings, with two
distinct realizations of optimality. For the $(3,4,4)$ octahedral
Coxeter tiling, the optimal packing density is less than the maximal
value, and there are three distinct configurations of balls yielding
the same optimal value. The case if the $(5,3,6)$ dodecahedral tiling we have obtained some interesting horoball arrangements with less densities.
Table 2 contains a summary of the optimal packing densities under our constraints.
\newpage
\begin{center}
Table 2. Optimal packing densities for the four fully asymptotic Coxeter
tilings
\vspace{3mm}

\begin{tabular}{|l|l|l}
\hline
Schl\"{a}fli symbol & Optimal density \\
\hline
$(3,3,6)$&  $0.853276^{*}$ \\
$(3,4,4)$&  $0.818808$ \\
$(4,3,6)$&  $0.853276^{*}$ \\
$(5,3,6)$&  $0.787251 $\\ 
\hline
\end{tabular}
\end{center}
\begin{center}
{$^{*}$These values are identical to the B\"or\"oczky and Florian limit.}
\end{center}

In the future it will be interesting to investigate tilings given
various uniform conditions on the configuration of the balls in
$\mathbb{H}^3$ as well as higher dimensional hyperbolic spaces.
These studies may show the existence of multiple optimal
configurations for given tilings, similarly as we have observed in
$\mathbb{H}^3$. To the knowledge of the authors the solution of the
above problem is still open.

Optimal sphere packings in other homogeneous Thurston geometries
represent another huge class of open mathematical problems. For
these non-Euclidean geometries only very few results are known
\cite{Sz07-2}, \cite{Sz10-1}. Detailed studies are the objective of
ongoing research.



\end{document}